\documentclass[12pt,reqno]{amsart}

\usepackage{diagrams,amssymb}

\diagramstyle[height=2em,width=2em, Postscript]
\newarrow{Epi}----{>>}
\newarrow{Mono}>--->
\newarrow{Iso}>---{>>}
\newarrow{Inc}C--->
\newarrow{Mapsto}|--->
\newarrow{Igual}=====
\newarrow{Dashto}{}{dash}{}{dash}>

\newcommand{\Z}{{\mathbb Z}}
\newcommand{\eS}{{\mathcal {S}}}

\newcommand{\F}{{\mathbb F}}

\newcommand{\oP}{{\overline{P}}}

\newcommand{\pcom}{_{p}^{\wedge}}

\newcommand{\B}[1]{B^{#1}\Z/p}
\newcommand{\bz}{B\Z/p}

\newcommand{\map}{\operatorname{map}\nolimits}

\newcommand{\A}{\ifmmode{\mathcal{A}}\else${\mathcal{A}}$\fi}
\newcommand{\K}{\ifmmode{\mathcal{K}}\else${\mathcal{K}}$\fi}
\newcommand{\U}{\ifmmode{\mathcal{U}}\else${\mathcal{U}}$\fi}
\newcommand{\T}{\ifmmode{\mathcal{T}}\else${\mathcal{T}}$\fi}

\newtheorem{Thm}{Theorem}[section]
\newtheorem{Prop}[Thm]{Proposition}
\newtheorem{Cor}[Thm]{Corollary}
\newtheorem{Lem}[Thm]{Lemma}

\theoremstyle{definition}

\newtheorem{Rmk}[Thm]{Remark}
\newtheorem{Ex}[Thm]{Example}

\theoremstyle{remark}

\title[Postnikov pieces and $B\Z/p$-homotopy theory]{Postnikov pieces and $B\Z/p
$-homotopy
theory}

\author{Nat\`{a}lia Castellana}
\author{Juan A. Crespo}
\author{J\'er\^{o}me Scherer}

\thanks{All three authors are partially supported by MEC grant MTM2004-06686, the third
author is supported by the program Ram\'on y Cajal, MEC, Spain.}

\begin{document}


\begin{abstract}
We present a constructive method to compute the cellularization
with respect to $\B{m}$ for any integer $m \geq 1$ of a large
class of $H$-spaces, namely all those which have a finite number
of non-trivial $\B{m}$-homotopy groups (the pointed mapping space
$\map_*( \B{m}, X)$ is a Postnikov piece). We prove in particular
that the $\B{m}$-cellularization of an $H$-space having a finite
number of $\B{m}$-homotopy groups is a $p$-torsion Postnikov
piece. Along the way we characterize the $B\Z/p^r$-cellular
classifying spaces of nilpotent groups.
\end{abstract}


\maketitle

\section*{Introduction}
\label{sec intro}

The notion of $A$-homotopy theory was introduced by Dror
Farjoun~\cite{Dror} for an arbitrary connected space $A$. Here $A$
and its suspensions play the role of the spheres in classical
homotopy theory and so the $A$-homotopy groups of a space $X$ are
defined to be the homotopy classes of pointed maps
$[\Sigma^iA,X]$. The analogue to weakly contractible spaces are
those spaces for which all $A$-homotopy groups are trivial. This
means that the pointed mapping space $\map_*(A,X)$ is
contractible, i.e. $X$ is an $A$-local space. On the other hand
the classical notion of $CW$-complex is replaced by the one of
$A$-cellular space. Such spaces that can be constructed from $A$
by means of pointed homotopy colimits.

Thanks to work of Bousfield~\cite{B2} and Dror Farjoun~\cite{Dror}
there is a functorial way to study $X$ through the eyes of $A$.
The nullification $P_AX$ is the biggest quotient of $X$ which is
$A$-local and $CW_AX$ is the best $A$-cellular approximation of
the space $X$. Roughly speaking $CW_A X$ contains all the
transcendent information of the mapping space $\map_*(A,X)$ since
it is equivalent to $\map_*(A,CW_A X)$. Hence explicit computation
of the cellularization would give access to information about
$\map_*(A,X)$. The importance of mapping spaces (in the case $A =
B\Z/p$) is well established since Miller's solution to the
Sullivan conjecture~\cite{Miller}.

While there is a lot of literature devoted to computations of $P_A
X$, only very few computations of $CW_A X$ are available. For
instance Chach\'olski describes a strategy to compute the cellularization $CW_{A}(X)$ in \cite{MR97i:55023}. This method has
been successfully applied in some cases (cellularization with
respect to Moore spaces \cite{MR2003a:55024},
$B\Z/p$-cellularization of classifying spaces of finite groups
\cite{Ramon}), but it is in general difficult to apply.

An alternative way to compute $CW_AX$ is the following. The
localization map $l:X\rightarrow P_AX$ provides an equivalence
$CW_AX\simeq CW_A\oP_A X$ where as usual $\oP _AX$ denotes the homotopy fiber
of $l$. This equivalence gives a strategy to compute
cellularizations when $CW_A\oP_AX$ is known. For instance if $X$
is $A$-local then $\oP _AX\simeq*$ and so $CW_AX\simeq *$. From
the $A$-homotopy point of view the next case in which the
$A$-cellularization should be accessible is when $X$ has only a
finite number of $A$-homotopy groups, that is some iterated loop
space $\Omega^nX$ is $A$-local.  Natural examples of spaces satisfying this condition
are obtained   by considering the $n$-connected covers of $A$-local spaces.

Let us specialize in $H$-spaces and $A=B^m\Z/p$. P.~Bousfield has
determined the fiber of the localization map $X\rightarrow
P_{\B{m}}X$ (see~\cite{B2}) when $\Omega ^n X$ is
$\B{m}$-local. He shows that for such an $H$-space $\oP _{\B{m}}X$
is a $p$-torsion Postnikov piece $F$, whose homotopy groups are
concentrated in degrees from $m$ to $m+n-1$. As $F$ is also an
$H$-space (because $l$ is an $H$-map), we call it an
\emph{H-Postnikov piece}. The cellularization of $X$ (which is
again an $H$-space because $CW_A$ preserves $H$-structures)
therefore coincides with that of a Postnikov piece. We do this in
Section~\ref{sec Postnikov} and this enables us to obtain our main
result.

\medskip

\noindent{\bf Theorem~\ref{cellularization of H-spaces}}
{\it Let $X$ be a connected $H$-space such that $\Omega ^nX$ is
$\B{m}$-local. Then $$CW_{\B{m}}X\simeq F\times K(W,m)$$ where $F$
is a p-torsion $H$-Postnikov piece with homotopy concentrated in
degrees from $m+1$ to $n$ and $W$ is an elementary abelian $p$-
group.}

\medskip

Thus, when $X$ is an $H$-space with only a finite number of
$\B{m}$-homotopy groups, the cellularization $CW_{\B{m}}X$ is a
$p$-torsion $H$-Postnikov piece! This is not true in general if we
do not assume $X$ to be an $H$-space. For instance, the
$\bz$-cellularization of $B\Sigma_3$ is a space with infinitely
many non-trivial homotopy groups~\cite{dichotomy}.

For $m=1$ there is a large class of $H$-spaces which is known to
have some local loop space by previous work of the authors
\cite{CCS}: those for which the mod $p$ cohomology is finitely
generated as an algebra over the Steenrod algebra. Hence we obtain
the following.

\medskip

\noindent{\bf Proposition~\ref{celularizacion con cohomologia-fg}}
{\it Let $X$ be a connected $H$-space such that $H^*(X;\F_p)$ is
finitely generated as algebra over the Steenrod algebra. Then
$$CW_{\bz}X\simeq F\times K(W,1)$$
where $F$ is a 1-connected $p$-torsion $H$-Postnikov piece and $W$
is an elementary abelian $p$-group. Moreover, there exists an integer $k$
such that $CW_{\B{m}}X\simeq *$ for $m\geq k$.}

\medskip

Our results allow explicit computations which we exemplify by
computing the $B\Z/p$-cellularization of the $n$-connected cover
of any finite $H$-spaces (Proposition~\ref{prop CWfiniteH}), as
well as the $\B{m}$-cellularizations of the classifying spaces for
real and complex vector bundles $BU$, $BO$, and their connected
covers $BSU$, $BSO$, $BSpin$, and $BString$, see
Proposition~\ref{prop BString}.


\section{A double filtration of the category of spaces}
\label{sec filtration}

As mentioned in the introduction the condition that $\Omega^nX$ be
$\B{m}$-local will enable us to compute the
$\B{m}$-cellurarization of $H$-spaces. This section is devoted to
give a picture of how such spaces are related for different
choices of $m$ and $n$.

First of all we present a technical lemma which collects various facts needed in the rest of the paper.

\begin{Lem}
\label{localization and loops}
 Let $X$ be a connected space and $m>0$. Then,
\begin{enumerate}
  \item  If $X$ is $\B{m}$-local then $\Omega^nX$ is
$\B{m}$-local for all $n\geq 1$.
  \item If $X$ is $\B{m}$-local  then it is $\B{m+s}$-local for all $s\geq 0$.
\item If $\Omega X$ is $\B{m}$-local, then
$X$ is $\B{m+s}$-local for all $s\geq 1$.
\end{enumerate}
\end{Lem}
\begin{proof} To prove (1) simply apply $\map_*(\bz,-)$ to the path
fibration $\Omega X \rightarrow *\rightarrow X$.

Statement (2) is given by Dwyer's version of Zabrodsky's lemma
\cite[Prop. 3.4]{MR97i:55028} to the universal fibration $\B{n}
\rightarrow *\rightarrow \B{n+1}$.

Finally (3) is a direct consequence of Zabrodsky's lemma (now in
its connected version~\cite[Prop. 3.5]{MR97i:55028}) applied to
the universal fibration and using the fact that $\Omega X$
$\B{n}$-local implies $\map(\B{n},X)_c\simeq X$.
\end{proof}

Of course the converses of the previous results are not true. For
the first statement take the classifying space of a discrete group
at $m=1$. For the second and third consider $X=BU$. It is a
$\B{2}$-local (see
Example~\ref{ExamplesofB(n+1)localwithloopsneverlocal}) space but
neither $BU$ not $\Omega BU$ are $\bz$-local. Observe that in fact
$\Omega ^nBU$ is never $\bz$-local. The next result shows that this is
the general situation for $H$-spaces. That is, if an $H$-space is
$\B{m+1}$-local then either $\Omega X$ is $\B{m}$-local or $\Omega
^nX$ is never $\B{m}$-local $\forall n\geq1$.

\begin{Thm}\label{LazosdeunB^n+1local}
Let $X$ be a $\B{n+1}$-local space such that $\Omega^kX$ is
$\B{n}$-local for some $k>0$. Then $\Omega X$ is $\B{n}$-local.
\end{Thm}

\begin{proof} It is enough to prove the result for $k=2$. Consider
the fibration
$$K(Q, n+1) \rTo P_{\Sigma^2 \B{n}} X \simeq X \rTo{} P_{\Sigma
\B{n}}X $$
where the fiber is a $p$-torsion Eilenberg-Mac Lane space by
Bousfield's description of the fiber of the $\Sigma
\B{n}$-nullification \cite[Theorem 7.2]{B2}. The total space is
$\B{n+1}$-local and so is the base by the previous lemma. Thus
$\map_*(\B{n+1}, K(Q, n+1))$ must be contractible as well, i.e. $Q
= 0$.
\end{proof}

The previous analysis leads to a double filtration of the category
of spaces. Let $n\geq 0$ and $m\geq 1$. We introduce the notation
$$
\mathcal S ^n_m=\{X;\; \Omega^nX \text{ is }
\B{m}\text{-local}\}\;.
$$

Lemma~\ref{localization and loops} yields then a diagram of
inclusions:

$$\begin{diagram}[Postscript]
  \eS^0_1 & \rInc & \eS^1_1 & \rInc & \eS^2_1& \rInc & \dots & \rInc & \eS^n_1 &
 \dots \\
  \dInc & \ldInc & \dInc &\ldInc  &\dInc&  &  &\ldInc& \dInc &  \\
  \eS^0_2 & \rInc & \eS^1_2 &\rInc & \eS^2_2  & \rInc & \dots & \rInc & \eS^n_2
& \dots \\
 \dInc & \ldInc & \dInc &\ldInc  &\dInc&  &  &\ldInc& \dInc &  \\
\vdots &  & \vdots &  & \vdots &  &  & &\vdots& \\
 \dInc & \ldInc & \dInc &\ldInc  &\dInc&  &  &\ldInc& \dInc &  \\
  \eS^0_m & \rInc & \eS^1_m &\rInc & \eS^1_m  & \rInc & \dots & \rInc &
  \eS^n_m
& \dots \\
\vdots &  & \vdots &  &\vdots  &  &  &  &\vdots&\\
\end{diagram}
$$

\begin{Ex}\label{CW y recubridores}
Examples of spaces in every stage of the filtration are known.

\begin{enumerate}
\item $\eS ^0_1$ are the spaces that are $\bz$-local.
This contains in particular any finite space (by Miller's theorem \cite[Thm.~A]{Miller}), and for a nilpotent
space $X$ (of finite type with finite fundamental group) to be $\bz$-local is equivalent to its cohomology
$H^*(X;\F_p)$ being locally finite by \cite[Corollary~8.6.2]{MR95d:55017}.
\item If $X\langle n \rangle $ denotes the $n$-connected
cover of a space $X$, then the homotopy fiber of $\Omega ^{n-1}
X\langle n\rangle \rightarrow \Omega ^{n-1} X$ is a discrete
space. Hence if $X\in \eS^0_m$ then $X\langle n\rangle \in \eS
^{n-1}_m$.
\item Observe that $\eS ^n_m \subset \eS^{n-k}_{m+k}$
for all $0\leq k\leq n$.
\item The previous examples provide
spaces in every stage of the double filtration. Consider a finite
space. It is automatically $\bz$-local. Its $n$-connected cover
$X\langle n\rangle$ lies in $\eS_1^{n-1}$. Hence $X\in
\eS_{k+1}^{n-k-1}$ for all $0\leq k\leq n$.
\end{enumerate}
\end{Ex}

Next example provides a number of spaces living in $\eS^0_{m}$ but
not obtained from the first row of the filtration by taking $n$-connected covers. Of course their
connected covers will be \emph{new} examples of spaces living in
$\eS^{n}_{m}$.

\begin{Ex}\label{ExamplesofB(n+1)localwithloopsneverlocal}
Let $E_*$ be a homology theory. If $\tilde{E}^i(K(\Z/p\Z,n))=0$
for $i\geq j$ then the spaces $E^i$ for $i\geq j$ representing the
corresponding homology theory are $\B{n}$-local. If
$\tilde{E}^j(K(\Z/p\Z,n-1))\neq 0$ then $E^j$ is not
$\B{n-1}$-local. In particular if $E_*$ is periodic, it follows
that the spaces $\{E^i\}$ for $i\geq j$ are $\B{n}$-local but none
of their iterated loops are $\B{n-1}$-local.

A first example of such a behavior is obtained from complex
K-theory, $BU$ is $\B{2}$-local but $BU$ and $U$ are not
$\bz$-local (see \cite{MR57:7585}). Note that real and
quaternionic $K$-theory enjoy the same properties.

For every $n$, examples of homology theories following this
pattern are given by $p$-torsion homology theories of type III-$n$
as described in \cite{MR57:17648}. The $n$th Morava $K$-theory
$K(n)_*$ for $p$ odd is an example of such behavior with respect to
Eilenberg-Mac~Lane spaces. The spaces representing $K(n)_*$ are
$\B{n+1}$-local but none of their iterated loops are
$\B{n}$-local.
\end{Ex}
Our aim is to provide tools to compute the $\B{m}$-cellularization
of any $H$-space lying in the $m$-th row of the above diagram. The
key point is the following result of Bousfield~\cite{B2}
determining the fiber of the localization map.

\begin{Prop}
\label{Bousfield2}
Let $n\geq 0$ and $X$ be a connected $H$-space $X$ such that
$\Omega^n X$ is $\B{m}$-local. Then there is an $H$-fibration
$$
F \rTo X \rTo P_{\B{m}} X
$$
where $F$ is a $p$-torsion $H$-Postnikov piece whose homotopy
groups are concentrated in degrees from $m$ to $m+n$. \hfill{\qed}
\end{Prop}

Therefore, since $F \rightarrow X$ is a $\B{m}$-cellular
equivalence, we only need to compute the cellularization of a
Postnikov piece (which will end up being a Postnikov piece again,
see Theorem~\ref{prop cell}). Actually even more is true.

\begin{Prop}
\label{prop converse}
Let $X$ be a connected space such that $CW_{\B{m}} X$ is a
Postnikov piece. Then there exists an integer $n$ such that
$\Omega^n X$ is $\B{m}$-local.
\end{Prop}

\begin{proof}
Let us loop once Chach\'olski fibration $CW_{\B{m}} X \rightarrow
X \rightarrow P_{\Sigma\B{m}} C$, see
\cite[Theorem~20.5]{MR97i:55023}. As $\Omega P_{\Sigma\B{m}} C$ is
equivalent to $P_{\B{m}} \Omega C$ by \cite[Theorem~3.A.1]{Dror},
we get a fibration over a $\B{m}$-local base space
$$
\Omega CW_{\B{m}} X \rTo \Omega X \rTo P_{\B{m}} \Omega C.
$$
Now there exists an integer $n$ such that $\Omega^n CW_{\B{m}} X$
is discrete, thus $\B{m}$-local. Therefore so is $\Omega^n X$.
\end{proof}

%
%

\section{Cellularization of fibrations over $BG$}
\label{sec fibrations}



In general the cellularization of the total space of a fibration
is very difficult to compute. We explain in this section how to
deal with this problem when the base space is the classifying
space of a discrete group. The first step applies to any group, in
the second, see Proposition~\ref{prop partialsections} below, we
specialize to nilpotent groups.

\begin{Prop}
\label{prop normal}
Let $r\geq 1$ and $F\rTo E \rTo^\pi BG$ be a fibration where $G$
is a discrete group. Let $S$ be the (normal) subgroup generated by
all elements $g\in G$ of order $p^i$ for some $i \leq r$ such that
the inclusion $B\langle g\rangle \rTo BG$ lifts to $E$ up to an
unpointed homotopy. Then the pullback of the fibration along $BS
\rTo BG$
$$
\begin{diagram}
E' & \rTo^f & E & \rTo^p & B(G/S) \\
\dTo & & \dTo_\pi & & \dTo\\
BS & \rTo & BG & \rTo^{p'} & B(G/S) \\
\end{diagram}
$$
induces a $B\Z/p^r$-cellular equivalence $f:E'\rightarrow E$ on
the total space level.
\end{Prop}

\begin{proof}
We have to show that $f$ induces a homotopy equivalence on pointed
mapping spaces $\map_*(B\Z/p^r,-)$. The top
fibration in the diagram yields a fibration
$$
\map_*(B\Z/p^r,E') \rTo^{f_*} \map_*(B\Z/p^r,E) \rTo^{p_*}
\map_*(B\Z/p^r,B(G/S)).
$$
Since the base is homotopically discrete we only need to check
that all components of the total space are sent by $p_*$ to the
component of the constant. Consider thus a map $h: B\Z/p^r
\rightarrow E$. The composite $p\circ h$ is homotopy equivalent to
a map induced by a group morphism $\alpha \colon \Z/p^r\rTo G$
whose image $\alpha(1)=g$ is in $S$ by construction. Therefore
$p\circ h = p'\circ \pi \circ h$ is
null-homotopic.
\end{proof}

\begin{Rmk}
If the fibration in the above proposition is an $H$-fibration (in particular
$G$ is abelian), the set of elements $g$ for which there is a lift to the
total space forms a subgroup of $G$. The central extension $Z(D_8) \hookrightarrow D_8
\rightarrow \Z/2 \times \Z/2$ of the dihedral group $D_8$ provides
an example where the subgroup $S$ is $\Z/2 \times \Z/2$ but
the element in $S$ represented by an element of order 4 in $D_8$
does not admit a lift.
\end{Rmk}

The next lemma is a variation of Dwyer's version of Zabrodsky's
Lemma in~\cite{MR97i:55028}.

\begin{Lem}\label{lemma kernel}
Let $F \rightarrow E \stackrel{f}{\longrightarrow} B$ be a
fibration over a connected base, and $A$ a connected space such
that $\Omega A$ is $F$-local. Then any map $g:E\rightarrow A$
which is homotopic to the constant when restricted to $F$ factors
through a map $h:B\rightarrow A$ up to unpointed homotopy and
moreover $g$ is pointed null-homotopic if and only if $h$ is so.
\end{Lem}

\begin{proof}
Since $\Omega A$ is $F$-local, we see that the component of the
constant map $\map_*(F, A)_c$ is contractible and therefore the
evaluation at the base point $\map(F, A)_c \rightarrow A$ is an
equivalence. By Proposition~3.5 in \cite{MR97i:55028}, $f$ induces
a homotopy equivalence
$$
\map(B,A)\simeq \map(E,A)_{[F]}.
$$
where $\map(E,A)_{[F]}$ denotes the space of those maps $E
\rightarrow A$ which are homotopic to the constant when restricted
to $F$.

We restrict now to the component of the constant map $c: E
\rightarrow A$. There is only one component in the pointed mapping
space sitting over $c$ since any map homotopic to the constant map
is also homotopic by a pointed homotopy. The result follows.
\end{proof}

\begin{Prop}
\label{prop partialsections}
Let $r \geq 1$ and $F\rTo^i E \rTo^\pi BG$ be a fibration where
$G$ is a nilpotent group generated by elements of order $p^i$ with
$i\leq r$. Assume that for each of these generators $x\in G$, the
inclusion $B\langle x\rangle \rTo BG$ lifts to $E$ up to unpointed
homotopy. If $F$ is $B\Z/p^r$-cellular then so is $E$.
\end{Prop}

\begin{proof}Chach\'olski's description \cite{MR97i:55023} of the
cellularization $CW_{B\Z/p^r}(E)$ as the homotopy fiber of the
composite $f: E \rTo C \rTo P_{\Sigma B\Z/p^r}(C)$ where $C$ is
the homotopy cofiber of the evaluation map $\bigvee_{[B\Z/p^r,E]}
B\Z/p^r\rTo E$ tells us that $E$ is cellular if the map $f$ is
null-homotopic. Observe that if $f$ is nullhomotopic then the fiber inclusion $CW_{B\Z/p^r}(E)\rTo E$ has a section and therefore $E$ is cellular since it is a retract of a cellular space (\cite[2.D.1.5]{Dror}).

As the existence of an unpointed homotopy to the constant map
implies the existence of a pointed one, we work now in the
category of unpointed spaces. Remark that for any map $g:Z\rTo E$
from a $B\Z/p^r$-cellular space $Z$, the composite $f\circ g$ is
null-homotopic since $g$ factors through the cellularization of
$E$. In particular the composite $f\circ i$ is null-homotopic. By
Lemma \ref{lemma kernel} there exists $\bar{f}:BG\rTo P_{\Sigma
B\Z/p^r}(C)$ such that $\bar{f}\circ \pi \simeq f$ and, moreover,
$f$ is null-homotopic if and only if $\bar{f}$ is so.


We first assume that $G$ is a finite group and show by induction
on the order of $G$ that $\bar f$ is null-homotopic. If $|G|=p$,
then the existence of a section $s:BG \rTo E$ implies that $f\circ
s=\bar{f}$ is nullhomotopic since $BG = B \Z/p$ is cellular.

Let $\{x_1,\ldots ,x_k\}$ be a minimal set of generators which
admit a lift. Let $H\trianglelefteq G$ be the normal subgroup
generated by $x_1,\ldots,x_{k-1}$ and their conjugates by powers
of~$x_k$. There is a short exact sequence $H\rTo G\rTo \Z/p^{a}$
where the quotient group is generated by the image of $x_k$.
Consider the fibration $F\rightarrow E'\rightarrow BH$ obtained by
pulling back along $BH\rightarrow BG$ and denote by
$h:E'\rightarrow E$ the induced map between the total spaces. The
inclusions in $G$ of two conjugate subgroups are (freely)
homotopic and so $H$ satisfies the assumptions of the proposition.
Thus the induction hypothesis tells us that $E'$ is cellular and
therefore $f\circ h$ is nullhomotopic. This implies that the
restriction of $\bar{f}$ to $BH$ is nullhomotopic. Consider the
following diagram
$$
\begin{diagram}
B(\langle x_k \rangle \cap H) & \rTo & BH &  & \\
\dTo & & \dTo & \rdTo^* & \\
B(\langle x_k\rangle) & \rTo & BG & \rTo^{\bar{f}} &  P_{\Sigma B\Z/p^r}(C)\\
\dTo & & \dTo & \ruTo_{f'}  & \\
B\Z/p^a & = & B\Z/p^a & &
\end{diagram}
$$
By Lemma \ref{lemma kernel}, it is enough to show that $f'$ is
nullhomotopic. But again by Lemma \ref{lemma kernel} applied to
the fibration on the left, we see that $f'$ is nullhomotopic since
$\bar{f}$ restricted to $\langle x_k\rangle$ is so. Therefore
$\bar{f}$ is nullhomotopic.

Assume now that $G$ is not finite. Any subgroup of $G$ generated
by a finite number of elements of order a power of $p$ has a
finite abelianization, and must therefore be finite itself by
\cite[Theorem~2.26]{MR48:11314}. Thus $G$ is locally finite, i.e.
$G$ is a filtered colimit of finite nilpotent groups generated by
elements of order $p^i$ for $i\leq r$. Likewise, $BG$ is a
filtered homotopy colimit of $BS$ where $S$ are finite groups
(generated by finite subsets of the set of generators) which
verify the hypothesis of the proposition. The total space $E$ can
be obtained as a pointed filtered colimit of the total spaces
obtained by pulling back the fibration. By the finite case
situation they are all cellular and therefore so is $E$.
\end{proof}

Sometimes the existence of the ``local" sections defined for
every generator permits to construct a global section of the fibration.
By a result of Chach\'olski \cite[Theorem~4.7]{MR97i:55023} the total space
of such a split fibration is cellular since $F$ and $BG$ are so.
This is the case for an $H$-fibration and $E$ is then weakly equivalent to
a product $F \times BG$.

A straightforward consequence of the above proposition (in the
case when the fibration is the identity of $BG$) is the following
characterization of the $B\Z/p^r$-cellular classifying spaces. For
$r=1$ we obtain R. Flores' result~\cite{Ramon}.

\begin{Cor}
\label{cor cellularBG}
Let $r \geq 1$ and $G$ be a nilpotent group generated by elements
of order $p^i$ with $i\leq r$. Then $BG$ is $B\Z/p^r$-cellular.
\hfill{\qed}
\end{Cor}

\begin{Ex}
\label{ex Q8}
The quaternion group $Q_8$ of order 8 is generated by elements of
order~4. Therefore $BQ_8$ is $B\Z/4$-cellular. We do not know an
explicit way to construct $BQ_8$ as a pointed homotopy colimit of
a diagram whose values are copies of $B\Z/4$.
\end{Ex}

The previous technical propositions allow us to state the main
result of this section. It provides a constructive description of
the cellularization of the total space of certain fibrations over
classifying spaces of nilpotent groups.

\begin{Thm}
\label{thm CWfibration}
Let $G$ be a nilpotent group and $F\rTo E\rTo BG$ be a fibration
with $B\Z/p^r$-cellular fiber $F$. Then the cellularization of $E$
is the total space of a fibration $F\rTo CW_{B\Z/p^r}(E)\rTo BS$
where $S \triangleleft G$ is the (normal) subgroup generated by
the $p$-torsion elements $g$ of order $p^i$ with $i\leq r$ such
that the inclusion $B\langle g \rangle \rTo BG$ lifts to $E$ up to
unpointed homotopy.
\end{Thm}

\begin{proof}
By Proposition~\ref{prop normal} pulling back along $BS
\rightarrow BG$ yields a cellular equivalence $f$ in the following
square:
$$
\begin{diagram}
E_S & \rTo ^f & E \\
\dTo & & \dTo \\
BS & \rTo & BG
\end{diagram}
$$
By Proposition~\ref{prop partialsections} the total space $E_S$ is
cellular and therefore $E_S \simeq CW_{B\Z/p^r}(E)$.
\end{proof}

\begin{Cor}
\label{cor finite}
Let $G$ be a nilpotent group and $S \triangleleft G$ be the
(normal) subgroup generated by the $p$-torsion elements $g$ of
order $p^i$ with $i\leq r$. Then $CW_{B\Z/p^r} BG \simeq BS$.
Moreover when $G$ is finitely generated, $S$ is a finite
$p$-group.
\end{Cor}

\begin{proof}
We only need to show that $S$ is a finite $p$-group. Notice that the abelianization of $S$ is $p$-torsion, then $S$ is also a torsion group (see \cite[Cor.~3.13]{Warfield}). Moreover, since $G$ is finitely generated, by \cite[3.10]{Warfield}, $S$ is finite.
\end{proof}

In fact in that case the previous result also holds when the base
space is an Eilenberg-Mac Lane space $K(G,n)$.

\begin{Prop}
Let $F\rTo^i E \rTo^\pi K(G,n)$ be a fibration where $G$ is a
finitely generated group by elements of order $p^i$ where $i\leq
r$ and $n>1$. Assume  that for each generator $x\in G$, the
inclusion $K(\langle x\rangle,n) \rTo K(G,n)$ lifts to $E$. If $F$
is $B\Z/p^r$-cellular then $E$ is so. \hfill{\qed}
\end{Prop}

\section{Cellularization of nilpotent Postnikov pieces}
\label{sec Postnikov}

In this section we compute the cellularization with respect to
$B\Z/p^r$ of nilpotent Postnikov pieces. The main difficulty lies
in the fundamental group, so it will be no surprise that these
results hold as well for cellularization with respect to
$B^m\Z/p^r$ with $m \geq 2$. We will often use the following
closure property \cite[Theorem~2.D.11]{Dror}.

\begin{Prop}
\label{prop closure}
Let $F \rightarrow E \rightarrow B$ be a fibration where $F$ and
$E$ are $A$-cellular. Then so is $B$. \hfill{\qed}
\end{Prop}

\begin{Ex}~\cite[Corollary~3.C.10]{Dror} The Eilenberg- Mac Lane
space $K(\Z/p^k, n)$ is $B\Z/p^r$-cellular for any integer $k$ and
any $n \geq 2$.
\end{Ex}

The construction of the cellularization is performed by looking
first at the universal cover of the Postnikov piece. We start with
the basic building blocs, the Eilenberg-Mac Lane spaces. For the
results on the structure on infinite abelian groups, we refer the
reader to Fuchs' book \cite{MR41:333}.

\begin{Lem}
\label{lemma EML}
An Eilenberg-Mac Lane space $K(A, m)$ with $m \geq 2$ is $B
\Z/p^r$-cellular if and only if $A$ is a $p$-torsion abelian
group.
\end{Lem}

\begin{proof}
That $A$ must be $p$-torsion is clear. Assume thus that $A$ is a
$p$-torsion group. If $A$ is bounded it is isomorphic to a direct
sum of cyclic groups. Since cellularization commutes with finite
products $K(A, m)$ is $B\Z/p^r$-cellular when $A$ is a finite
direct sum of cyclic groups. Taking a (possibly transfinite)
telescope of $B\Z/p^r$-cellular spaces we obtain that $K(A, m)$ is
so for any bounded group.

In general $A$ splits as a direct sum of a divisible group $D$ and
a reduced one $T$. A $p$-torsion divisible group is a direct sum
of copies of $\Z/p^\infty$, which is a union of bounded groups,
thus $K(D, m)$ is cellular. Now $T$ has a basic subgroup $P < T$
which is a direct sum of cyclic groups and the quotient $T/P$ is
divisible. So $K(T, m)$ is the total space of a fibration
$$
K(P, m) \rightarrow K(T, m) \rightarrow K(D, m)
$$
When $m \geq 3$ we are done because of the above mentioned closure
property Proposition~\ref{prop closure}. If $m = 2$ we have to
refine the analysis of the fibration because $K(D, m-1)$ is not
cellular. However, as $D$ is a union of bounded groups $D[p^k]$,
the space $K(T, 2)$ is the telescope of total spaces $X_k$ of
fibrations with cellular fiber $K(P, 2)$ and base $K(D[p^k], 2)$.
We claim that these total spaces are cellular (and thus so is
$K(T, 2)$) and proceed by induction on the bound. Consider the
subgroup $D[p^{k}] < D[p^{k+1}]$ whose quotient is a direct sum of
cyclic groups $\Z/p$. Therefore $X_{k+1}$ sits in a fibration
$$
K(\oplus \Z/p, 1) \rightarrow X_{k} \rightarrow X_{k+1}
$$
where fiber and total space are cellular. We are done.
\end{proof}

We are now ready to prove that any $p$-torsion simply connected
Postnikov piece is a $B\Z/p^r$-cellular space.

\begin{Prop}
\label{proposition 1-connected}
A simply connected Postnikov piece is $B\Z/p^r$-cellular if and
only if it is $p$-torsion.
\end{Prop}

\begin{proof}
Let $X$ be a simply connected $p$-torsion Postnikov piece. For
some integer $m$, the $m$-connected cover $X\langle m \rangle$ is
an Eilenberg-Mac Lane space, which is cellular by Lemma~\ref{lemma
EML}. Consider the principal fibration
$$
K(\pi_m(X),m-1) \rTo X\langle m \rangle \rTo X\langle m-1 \rangle
$$
If $m \geq 3$ both $X\langle m\rangle$ and $K(\pi_m(X),m-1)$ are
cellular. It follows that $X\langle m-1 \rangle $ is cellular by
the closure property Proposition~\ref{prop closure}. The same
argument shows that $X\langle 2 \rangle$ is cellular.

Let us thus look at the fibration $X\langle 2 \rangle \rightarrow
X \rightarrow K(\pi_{2} X, 2)$. The same discussion on the
$p$-torsion group $\pi_{2} X$ as in Lemma~\ref{lemma EML} will
apply. If this is a bounded group, say by $p^k$, an induction on
the bound shows that $X$ is actually the base space of a fibration
where the total space is cellular because its second homotopy
group is $p^{k-1}$-bounded, and the fiber is cellular because it
is of the form $K(V, 1)$ with $V$ a $p$-torsion abelian groups
whose torsion is bounded by $p^r$. The closure property
Proposition~\ref{prop closure} ensures that $X$ is then cellular.

If $\pi_{2} X$ is divisible, $X$ is a telescope of cellular
spaces, hence cellular. If it is reduced, taking a basic subgroup
$B < \pi_{2} X$ yields a diagram of fibrations
$$
\begin{diagram}
X\langle 2 \rangle &\rTo&Y&\rTo &K(B, 2)\\
\dTo&&\dTo&&\dTo\\
X\langle 2 \rangle&\rTo&X&\rTo&K(\pi_{2} X, 2)\\
\dTo&&\dTo&&\dTo\\
*&\rTo&K(D, 2)&\rIgual &K(D, 2)
\end{diagram}
$$
which exhibits $X$ as the total space of a fibration over $K(D,
2)$ with $D$ divisible and $B\Z/p^r$-cellular fiber. Therefore
writing $D$ as a union of bounded groups as in the proof of Lemma
\ref{lemma EML}, $X$ is a telescope of cellular spaces, therefore
it is $B\Z/p^r$-cellular as well.
\end{proof}

\begin{Rmk}
The proof of the proposition holds in the more general setting
when $X$ is a $p$-torsion space such that $X\langle m \rangle$ is
$B\Z/p^r$ cellular for some $m\geq 2$. The proposition corresponds
to the case when some $m$-connected cover $X\langle m \rangle$ is
contractible.
\end{Rmk}

Recall from \cite[Corollary 2.12]{MR57:17635} that a connected
space is nilpotent if and only if its Postnikov system admits a
principal refinement $\dots \rightarrow X_s \rightarrow X_{s-1}
\rightarrow \dots \rightarrow X_1 \rightarrow X_0$. This means
that each map $X_{s+1} \rightarrow X_{s}$ in the tower is a
principal fibration with fiber $K(A_s, i_s-1)$ for some increasing
sequence of integers $i_s \geq 2$. We are only interested in
finite Postnikov pieces, i.e. nilpotent spaces that can be
constructed in a finite number of steps by taking homotopy fibers
of $k$-invariants $X_{s} \rightarrow K(A_s, i_s)$.

The key step in the study of the cellularization of a nilpotent
finite Postnikov piece is the analysis of a principal fibration
(given in our case by the $k$-invariants).

\begin{Thm}\label{prop cell}
Let $X$ be a $p$-torsion nilpotent Postnikov piece. Then there
exists a fibration
$$
X\langle 1 \rangle \rTo CW_{B\Z/p^r}X\rTo BS
$$
where $S$ is the (normal) subgroup of $\pi_1(X)$ generated by the
elements $g$ of order $p^i$ with $i \leq r$ such that the
inclusion $B \langle g \rangle \rightarrow B\pi_1 X$ admits a lift
to $X$ up to unpointed homotopy.
\end{Thm}

\begin{proof}
By Proposition \ref{proposition 1-connected} the universal cover
$X\langle 1 \rangle$ is cellular and there is a fibration $X
\langle 1 \rangle\rTo X\rTo BG$ where $G=\pi_1(X)$ is nilpotent.
The result follows then from Theorem~\ref{thm CWfibration}.
\end{proof}

\section{Cellularization of $H$-spaces}
\label{sec Cel of H-spaces}
In this section we will apply the computations of the
cellularization of $p$-torsion nilpotent Postnikov systems to
determine $CW_{\B{}}X$ when $X$ is an $H$-space. We prove:

\begin{Thm}
\label{cellularization of H-spaces}
Let $X$ be a connected $H$-space such that $\Omega^nX$ is
$\B{}$-local. Then
$$
CW_{\B{}}X\simeq Y\times K(W,1)
$$
where $Y$ is a simply connected p-torsion $H$-Postnikov piece with
homotopy concentrated in degrees $\leq n$ and $W$ is an elementary
abelian $p$-group.
\end{Thm}

\begin{proof}
The fibration in Bousfield's result Proposition~\ref{Bousfield2}
yields a cellular equivalence between a connected $p$-torsion
$H$-Postnikov piece $F$ and $X$. Theorem~\ref{prop cell} thus
applies. Moreover, as $F$ is an $H$-space as well, the subgroup
$S$ is abelian generated by elements of order $p$. Therefore we
have a fibration $F\langle 1 \rangle \rightarrow
CW_{\bz}F\rightarrow K(W,1)$ which admits a section (summing up the
local section). The cellularization therefore splits.
\end{proof}

This result applies for $H$-spaces satisfying certain finiteness
conditions.

\begin{Prop}
\label{celularizacion con cohomologia-fg}
Let $X$ be a connected $H$-space such that $H^*(X;\F_p)$ is
finitely generated as algebra over the Steenrod algebra. Then
$$
CW_{\bz}X\simeq F\times K(W,1)
$$
where $F$ is a 1-connected $p$-torsion $H$-Postnikov piece and $W$
is an elementary abelian $p$-group. Moreover there exists an integer $k$
such that $CW_{\B{m}}X\simeq *$ for $m\geq k$.
\end{Prop}

\begin{proof}
In~\cite{CCS} the authors prove that if $H^*(X;\F _p)$ is finitely
generated as algebra over the Steenrod algebra then $\Omega^nX$ is
$\bz$-local for some $n\geq 0$. Hence Theorem~\ref{cellularization
of H-spaces} applies and we obtain the desired result. In addition
Lemma~\ref{localization and loops} shows that $X$ is
$\B{n+s+1}$-local for any $s\geq 0$, which implies the second part
of the result.
\end{proof}

The technique we propose in this paper is not only a nice
theoretical tool which provides a general statement about how the
$B\Z/p$-cellularization of $H$-spaces look like. Our next result
shows that one can actually identify precisely this new space when
dealing with connected covers of finite $H$-spaces. Recall that by
Miller's theorem \cite[Thm.~A]{Miller}, any finite $H$-space $X$
is $\bz$-local and therefore $CW_{\bz}(X) \simeq *$. The universal
cover of $X$ is still finite and thus $CW_{\bz}(X\langle 1
\rangle)$ is contractible as well. We can therefore assume that
$X$ is $1$-connected. The computation of the cellularization of
the $3$-connected cover is already implicit in \cite{BC}.

\begin{Prop}
\label{prop CWfiniteH}
Let $X$ be a simply connected finite $H$-space and let $k$ denote
the rank of the free abelian group $\pi_3 X$. Then $CW_{B\Z/p} X
\langle 3 \rangle \simeq K(\oplus_k \Z/p, 1)$. For $n \geq 4$, up
to $p$-completion, the universal cover of $CW_{\bz}(X\langle n
\rangle)$ is weakly equivalent to the $2$-connected cover of
$\Omega (X[n])$.
\end{Prop}

\begin{proof}
By Browder's famous result \cite[Theorem~6.11]{MR23:A2201} $X$ is
even $2$-connected and its third homotopy group $\pi_3 X$ is free
abelian (of rank $k$) by J.~Hubbuck and R.~Kane's theorem
\cite{MR53:1582}. This means we have a fibration
$$
K(\oplus_k \Z_{p^\infty}, 1) \rTo X \langle 3 \rangle \rTo
P_{B\Z/p} X \langle 3 \rangle
$$
which shows that $CW_{B\Z/p} X \langle 3 \rangle \simeq K(\oplus_k
\Z/p, 1)$. We deal now with the higher connected covers. Consider
the following commutative diagram of fibrations
$$
\begin{diagram}
F & \rTo & F & \rTo & * \\
\dTo & & \dTo & & \dTo \\
\Omega X[n] & \rTo & X\langle n \rangle  & \rTo & X \\
\dTo & & \dTo & & \dIgual \\
P_{\bz}(\Omega X[n]) & \rTo & P_{\bz}(X\langle n \rangle) & \rTo & X \\
\end{diagram}
$$
where $F$ is a $p$-torsion Postnikov piece by \cite[Thm 7.2]{B2}
and the fiber inclusions are all $\bz$-cellular equivalences
because the base spaces are $B\Z/p$-local. Therefore
$$
CW_{\bz}(X\langle n \rangle)\simeq CW_{\bz}(F)\simeq F\langle 1
\rangle \times K(W,1)
$$
We wish to identify $F\langle 1 \rangle$. Since the fibrations in
the diagram are nilpotent, by \cite[II.4.8]{MR51:1825} they remain
fibrations after $p$-completion. By Neisendorfer's theorem
\cite{MR96a:55019} the map $P_{\bz}(X\langle n \rangle)
\rightarrow X$ is an equivalence up to $p$-completion, which means
that $P_{\bz}(\Omega (X[n]))\pcom \simeq *$. Thus $F\pcom \simeq
(\Omega (X[n]))\pcom$.  Notice that $\Omega (X[n])$ is
$1$-connected and its second homotopy group is free by the above
mentioned theorem of Hubbuck and Kane (which corresponds up to
$p$-completion to the direct sum of $k$ copies of the Pr\"ufer
group $\Z/p^\infty$ in $\pi_1 F$). Hence $F\langle 1 \rangle$
coincides with $\Omega (X[n])\langle 2\rangle$ up to
$p$-completion.
\end{proof}

To illustrate this result we compute the $B\Z/2$-cellularization
of the successive connected covers of $S^3$. The only delicate
point is the identification of the fundamental group.

\begin{Ex}
\label{sphere}
Recall that $S^3$ is $B\Z/2$-local since it is a finite space.
Thus the cellularization $CW_{B\Z/2} S^3$ is contractible. Next the fibration
$$
K(\Z_{2^\infty}, 1) \rightarrow S^3\langle 3 \rangle \rightarrow
P_{B\Z/2} S^3\langle 3 \rangle
$$
shows that $CW_{B\Z/2} S^3\langle 3 \rangle \simeq K(\Z/2, 1)$.
Finally since $S^3[4]$ does not split as a product (the
$k$-invariant is not trivial), we see that $CW_{B\Z/2} S^3\langle
4 \rangle \simeq K(\Z/2, 3)$. Likewise, for any integer $n \geq
4$, we have that $CW_{B\Z/2} S^3\langle n \rangle$ is equivalent
to the $2$-completion of the $2$-connected cover of $\Omega
(S^3[n])$. The same phenomenon occurs at odd primes.
\end{Ex}

\section{Cellularization with respect to $\B{m}$}
\label{sec higher}

All the techniques developed for fibrations over $BG$ apply to
fibrations over $K(G,n)$ when $n>1$ and we get the following
results.

\begin{Lem}
Let $n \geq 2$ and $X$ be a connected space. Then
$$CW_{B^n\Z/p^r}(X)=CW_{B^n\Z/p^r}(X\langle n-1 \rangle)\;.$$
\end{Lem}

\begin{proof}
For $n > i$ consider the following fibrations over
$B^n\Z/p^r$-local base spaces:
$$
X\langle i \rangle \rTo X\langle i-1 \rangle \rTo K(\pi_{i}(X),i)
$$
We see that $CW_{B^n\Z/p^r}(X \langle i \rangle)=CW_{B^n\Z/p^r}(X
\langle i-1 \rangle)$.
\end{proof}

\begin{Prop}
\label{prop cell2}
Let $m \geq 2$ and $X$ be a $p$-torsion nilpotent Postnikov piece.
Then there exists a fibration
$$
X\langle m \rangle\rTo CW_{\B{m}} X \rTo K(W,m)
$$
 where $W$ is a $p$-torsion subgroup of $\pi_m(X)$
with torsion bounded by~$p^r$.
\end{Prop}

\begin{Thm}
\label{cellularization of H-spaces 2}
Let $X$ be a connected $H$-space such that $\Omega ^nX$ is
$\B{m}$-local. Then
$$
CW_{\B{m}}X\simeq F\times K(W,m)
$$
where $F$ is a p-torsion $H$-Postnikov piece with homotopy
concentrated in degrees from $m+1$ to $n$ and $W$ is an elementary
abelian $p$-group.
\end{Thm}

\begin{Ex}
\label{Milgram}
Let $X$ denote ``Milgram's space", see \cite{MR40:2061}, the fiber
of $Sq^2: K(\Z/2, 2) \rightarrow K(\Z/2, 4)$. This is an infinite
loop space. By Proposition~\ref{proposition 1-connected} we know
it is already $B\Z/2$-cellular. Let us compute the cellularization
with respect to $B^m\Z/2$ for higher $m$'s. Since the
$k$-invariant is not trivial, we see that $CW_{B^2\Z/2} X \simeq
CW_{B^3\Z/2} X \simeq K(\Z/2, 3)$.
\end{Ex}

We compute finally the cellularization of the (infinite loop)
space $BU$ and its $2$-connected cover $BSU$ with respect to
Eilenberg-MacLane spaces $\B{m}$. By Bott periodicity this
actually tells us the answer for all connected covers of $BU$.

\begin{Ex}
\label{BU}
First of all, recall from Example
\ref{ExamplesofB(n+1)localwithloopsneverlocal} that $BU$ is
$\B{2}$-local since $\widetilde{K}^*(\B{2})=0$ and its iterated
loops are never $B\Z/p$-local. Therefore $CW_{\B{m}}(BU)$ is
contractible if $m\geq 2$. Since $BU\simeq BSU\times BS^1$ the
same property holds for $BSU$.
\end{Ex}

We now compute the $\B{m}$-cellularization of $BO$ and its
connected covers $BSO$, $BSpin$, and $BString$.

\begin{Prop}
\label{prop BString}
Let $m \geq 2$. Then
\begin{itemize}
\item[(i)] $CW_{\B{m}}(BO) \simeq CW_{\B{m}}(BSO) \simeq CW_{\B{m}}(BSpin) \simeq*$,
\item[(ii)] $CW_{\B{m}}(BString) \simeq *$ if $m > 2$,
\item[(iii)] $CW_{\B{2}}(BString) \simeq K(\Z/p, 2)$ and $\map_*(\B{2}, BString) \simeq
\Z/p$.
\end{itemize}
\end{Prop}

\begin{proof}
In \cite{MR80b:55003} W. Meier proves that real and complex
$K$-theory have the same acyclic spaces, hence $BO$ is also
$\B{2}$-local. Therefore $CW_{\B{m}}(BO)$ is contractible for any
$m\geq 2$. The $2$-connected cover of $BO$ is $BSO$ and there is a
splitting $BO\simeq BSO\times B\Z/2$, so that $CW_{\B{m}}(BSO)
\simeq *$.

The $4$-conected cover of $BO$ is $BSpin$. It follows from the
fibration
$$
BSpin\rightarrow BSO \stackrel{w_2}{\rightarrow} K(\Z/2,2)
$$
that that fiber of $BSpin\rightarrow BSO$ is $B\Z/2$. Since $BSO$
and $B\Z/2$ are $\B{2}$-local, so is $BSpin$. Therefore
$CW_{\B{m}}(BSpin)$ is contractible.

Finally, the $8$-connected cover of $BO$ is $BString$. It is the
homotopy fiber of $BSpin\stackrel{p_1/4}{\rightarrow} K(\Z,4)$,
where $p_1$ denotes the first Pontrjagin class. Consider the
fibration
$$
K(\Z,3)\rightarrow BString\rightarrow BSpin
$$
where the base space is $\B{m}$-local for $m\geq 2$. Together with
the exact sequence $\Z\rightarrow \Z[\dfrac{1}{p}]\rightarrow
\Z/p^\infty$, this implies that
$$
CW_{\B{m}}(BString) \simeq CW_{\B{m}}(K(\Z,3)) \simeq
CW_{\B{m}}(K(\Z/p^\infty,2))
$$
which is contractible unless $m=2$, when we obtain $K(\Z/p, 2)$.
This yields the explicit description of the pointed mapping space
$\map_*(\B{2}, BString)$.
\end{proof}

Observe that the iterated loops of the $m$-connected covers of
$BO$ and $BU$ are never $B\Z/p$-local. Hence we know that their
cellularization with respect to $B\Z/p$ must have infinitely many
non-vanishing homotopy groups by Proposition~\ref{prop converse}.


\bibliographystyle{amsplain}
\bibliography{ccs}



\bigskip
{\small
\begin{center}
\begin{minipage}[t]{8 cm}
Nat\`{a}lia Castellana and J\'er\^{o}me Scherer\\ Departament de
Matem\`atiques,\\ Universitat Aut\`onoma de Barcelona,\\ E-08193
Bellaterra, Spain
\\\textit{E-mail:}\texttt{\,natalia@mat.uab.es}, \\
\phantom{\textit{E-mail:}}\texttt{\,\,jscherer@mat.uab.es}
\end{minipage}
\begin{minipage}[t]{7 cm}
Juan A. Crespo \\ Departament de Economia i de Hist\`{o}ria
Econ\`{o}mica,
\\ Universitat Aut\`onoma de Barcelona,\\ E-08193 Bellaterra,
Spain
\\\textit{E-mail:}\texttt{\,JuanAlfonso.Crespo@uab.es},
\end{minipage}
\end{center}}

\end{document}